\documentclass{kapproc}
\usepackage{amssymb,amsmath}
\usepackage{array}
\usepackage{rotating}

\newcommand{\beq}{\begin{equation}}
\newcommand{\eeq}{\end{equation}}

\begin{document}
%\inputencoding{cp1251}

\normallatexbib

%\begin{center}
%\vspace{1cm}
%{\Large {\bf
\articletitle{Noncompact quantum\\
algebra \mbox{\lowercase{$u_q(2,1)$}}: \\
Positive discrete series\\ of irreducible representations }
%\vspace{1cm}
%\noindent
 \chaptitlerunninghead{Noncompact quantum algebra $u_q(2,1)$}
\author{Yu. F. Smirnov}
\affil{Skobeltsyn Institute of Nuclear Physics,
Moscow State University,
Moscow, 119992, Russia}
%and
\author{Yu. I. Kharitonov\thanks{Deceased.}}
%\footnote{\it Petersburg Institute of Nuclear Physics, Gatchina.} %}

% Nuclear Physics Institute of Moscow State University
%\end{center}

%\noindent
%\vspace*{0.5cm}

%{\bf A b s t r a c t}

%\vspace*{0.3cm}
\begin{abstract}
The structure positive of unitary irreducible
representations of the noncompact $u_q(2,1)$ quantum algebra
that are related to a positive discrete series is examined.
With the aid of projection operators for the $su_q(2)$ subalgebra,
a $q$-analog of the Gel'fand--Graev formulas is derived in the
basis corresponding to the reduction
$u_q(2,1)\rightarrow su_q(2)\times u(1)$.
Projection operators for the $su_q(1,1)$ subalgebra are employed to study
the same representations for the reduction $u_q(2,1)\rightarrow u(1)\times su_q(1,1)$.
The matrix elements of the generators of the $u_q(2,1)$ algebra
are computed in this new basis. A general analytic expression
for an element of the transformation bracket $\langle U|T\rangle_q$ between the
bases associated with above two reductions (the elements of this matrix are referred
to as $q$-Weyl coefficients) is obtained for a general case where
the deformation parameter $q$ is not equal to a root of unity.
It is shown explicitly that, apart from a phase, $q$-Weyl coefficients
coincide with the $q$-Racah coefficients for the $su_q(2)$ quantum algebra.
\end{abstract}

%\vspace*{0.5cm}

\section{Introduction}

  It is well known that the group theory methods are widely used in
  the theory of nucleus. They form the basis of nuclear spectroscopy
  and of various nuclear models, including the shell model, models
  dealing with collective degrees of freedom, and the interacting boson model.
   Since group theory or algebraic models usually admit an analytic
  solution, they are employed to study various properties of nuclear systems
  in particular, some of their asymptotic properties.
 For example, the popular Elliott model, based on $SU(3)$ symmetry, was
 successfully employed by Belyaev and his colleagues~\cite{1s} to analyze the
 asymptotic properties of the generalized density matrix.
   The discovery of quantum algebras and groups that was made more than 20
 years  by mathematical physicists of Leningrad school \cite{2s}--\cite{4s}
 gave new impetus to the development of algebraic methods in
  theoretical physics,
 in particular, to searches for applications of the representations of quantum groups
 and algebras in physics. For example, the construction of $q$ analogs of various
 nuclear models became a new field of research in theoretical nuclear
  physics.
  The point is that quantum algebras involve an additional variable parameter,
  the deformation parameter $q$. This renders models based
   on quantum algebras more adaptable and extends their potential in describing
physical systems (see, for example, the study of Raychev et al~\cite{5s}
and the review article of Bonatsos and Daskaloyannis  ~\cite{6s})
which is devoted to applications of quantum algebras in the theoretical nuclear
physics).
   However, the searches for physical applications of quantum
algebras must be preceded by a detailed investigations of their irreducible
representations.
 In this connection, the structure of unitary irreducible representations of the compact
$u_q ( 3 )$ algebra was examined in details in ~\cite{1}--\cite{10}. In our opinion
it is important to extend these results on the noncompact
$u_q(2,1)$ quantum algebra. The classical algebra $u(2, 1)$ describes
the dynamical symmetry
of a two-dimensional harmonic oscillator and of some other physical systems.
In view of this, a comprehensive analysis of unitary irreducible representations
of its quantum analogs may be helpful in constructing respective physical models.
In the present study, we restrict our consideration to the unitary irreducible
representations associated with a positive discrete series.

 Unitary irreducible representations of conventional (nondeformed)\linebreak
$u(n,m)$ algebras were studied by Gel'fand and Graev \cite{11} (see also
\cite{12}), who showed, among other things, that the unitary
irreducible representations of the $u(2,1)$ algebra can be divided into
three discrete series. The series of the highest weight unitary irreducible
representations or a negative discrete series consists of
representations such that each includes the highest weight vector
$|H\rangle$ that is, a vector annihilated by any raising generator of
the algebra.

  A positive discrete series is the series of representations such
that each includes the lowest weight vector $|L\rangle$ that is, a vector annihilated
by any lowering generator. There is yet another series, that is referred to as
an intermediate one and which is formed by unitary irreducible representations
having neither the highest weight vector $|H\rangle$ nor the lowest weight vector $|L\rangle$.
For this reason, this series deserves a dedicated consideration.

  Gel'fand and Graev \cite{11} presented explicit expressions for the matrix elements
 of generators associated with the above representations that is, the matrix elements
 of the generators $A_{ik}$ of the $u(n,m)$ algebra in the basis corresponding to
 the following reduction of this algebra to the chain of subalgebras:
 \beq
u(n,m)\rightarrow u(n,m-1)\rightarrow \ldots \rightarrow
u(n)\rightarrow \ldots \rightarrow u(2)\rightarrow u(1),
n\geq m.
\eeq
  However these authors did not give a regular procedure for deriving
  the expressions that they quoted in \cite{11}. For the $u(n,1)$,
these formulas were derived in \cite{13}--\cite{15}, but there is no
derivation of such formulas for the general case of the $u(n,m)$
algebras.
 In this study, we extend,  the approach proposed by Vilenkin in
 \cite{16} for the case of the $u(2,1)$ classical algebra  and
 examine the structure of its unitary irreducible representations
 associated with the positive series.
  These results obtained in this way are readily generalized
  to the case of negative discrete series. The intermediate discrete
  series of unitary irreducible representations will be
 considered in a separate paper.
 As in \cite{1}--\cite{10}, we assume
 that the deformation parameter $q$ is specified by an arbitrary
 positive number and define $q$-numbers and $q$-factorials
  as follows:
\beq
[n]=\frac{q^n-q^{-n}}{q-q^{-1}}\,,
\eeq
\beq
[n]!=[n][n-1] \ldots [2][1]\,, \quad [0]!=1.
\eeq
Below, we employ brackets to denote $q$-numbers, enclose
the signatures of unitary irreducible representations in Dirac
brackets, and reserve parentheses for the weight of a vector, for
example, the symbol $|\langle f_1f_2f_3\rangle (m_1m_2m_3)\rangle$
 stands for a
basis vector of a weight  $(m_1m_2m_3)$ in the unitary irreducible
representation $D^{\langle f_1f_2f_3\rangle }=D^{\langle f\rangle }$.

 \section{Positive discrete series of unitary irreducible representations}

   The $u(2,1)$ algebra is known to involve nine generators $A_{ik}$
$(i,k=1,2,3)$ satisfying the same commutation relations that the corresponding generators
of the compact $u(3)$ classical Lie algebra. However, properties of the $u(2,1)$
generators under Hermitian conjugations differ from those of the
$u(3)$ generators. The ``compact'' generators  $A_{11}$,
$A_{22}$, $A_{33}$, $A_{12}$ and $A_{21}$ of the $u(2,1)$ algebra have the same Hermitian
properties, as the $u(3)$ generators,
\beq
A_{ik}^+=A_{ki}\,,
\eeq
whereas the ``noncompact'' generators $A_{13}$, $A_{23}$, $A_{31}$ and
$A_{32}$ satisfy the relations
\beq
A_{13}^+=-A_{31}\,,
\eeq
\beq
A_{23}^+=-A_{32}\,.
\eeq
The minus sign in formulas (5) and (6) generates a fundamental distinction
between the structure of any unitary irreducible representation of the $u(2,1)$
algebra and the structure of the corresponding unitary irreducible representation
of the $u(3)$ algebra: all unitary irreducible representations of the compact
$u(3)$ algebra are finite-dimensional, whereas all unitary irreducible representations
of the noncompact $u(2,1)$ algebra(with the exception of the trivial identity
representation ) are infinite-dimensional.
  The noncompact  $u_q(2,1)$ quantum algebra is also specified by
nine generators $A_{ik}$ $(i,k=1,2,3)$ satisfying the same
commutation relations as the generators of the $u_q(3)$ compact quantum algebra.
The explicit expressions for these commutators can be found in \cite{1}.

 As to their properties with respect to Hermitian conjugation, those in (4) and (6)
 remain valid, whereas, in view of the relations
\beq
A_{13}^+=\tilde A_{31}=A_{32}A_{21}-qA_{21}A_{32}\neq A_{31}\,,
\eeq
\beq
A_{31}^+=\tilde
A_{13}=A_{12}A_{23}-q^{-1}A_{23}A_{12}\neq A_{13}\,,
\eeq
for the $u_q(3)$ algebra, that in (5) must be replaced by
\beq
A_{13}^+=-\tilde A_{31}\,.
\eeq

With the aid of (7) and (8), this relation can be recast into
either of the following two equivalent forms:
\beq
A_{13}^+=-A_{31}+(q-q^{-1})A_{21}A_{12}\,,
\eeq
or
\beq
A_{13}^+=-q^2A_{31}+(1-q^2)A_{32}A_{21}\,.
\eeq

For the  $u_q(2,1)$ algebra, we will consider the unitary irreducible
representation $D^{\langle f\rangle }$ of the lowest weight $(f)=(f_1f_2f_3)$:that is,
we assume that, in the space of this representation, there is the lowest
weight vector $|L\rangle$ that satisfies the relations
\beq
A_{ii}|L\rangle =f_i|L\rangle \,,\qquad (i=1,2,3)
\eeq
 annihilated by a pair of lowering generators
\beq
A_{31}|L\rangle =0\quad\mbox{and}\quad A_{32}|L\rangle =0 \,.
\eeq
Also it is annihilated by one raising generator
\beq
A_{12} |L\rangle =0.
\eeq

It is assumed that this vector is normalized by the relation
\beq
\langle L|L\rangle=1\,.
\eeq
All the other basis vectors  $|X\rangle$ of this unitary irreducible
representation  can be derived by applying the generators
$A_{13},A_{23}$ and $A_{21}$ to $|L\rangle $
\beq
|X\rangle =A_{21}^gA_{23}^kA_{13}^\ell|H\rangle .
\eeq

  In order to construct a basis of any unitary irreducible
representation of the $u_q(2,1)$ algebra, it is necessary to
specify a chain of subalgebras, and this can be done,
as it is well known, in three ways. The first way is to use
the $U$-spin subalgebra involving the generators $A_{11},
A_{12}, A_{21}$, and $A_{22}$, in which case the respective
reduction is
\beq
u_q(2,1)\rightarrow u_q(2)\rightarrow u_q(1)\,.
\eeq
One can also use the generators $A_{22}$, $A_{23}$, $A_{32}$, and
 $A_{33}$ forming the basis of the $T$-spin subalgebra or
the generators $A_{11}, A_{13}, A_{31}$, and $A_{33}$ generating
the $V$-spin subalgebra. Either of these two subalgebras
correspond to the reduction
\beq
u_q(2,1)\rightarrow u_q(1,1)\rightarrow u_q(1)\,.
\eeq
In this study, we restrict our consideration to the case  of $U$ and
$T$-spin bases.

\section{Basis vectors and matrix elements of the generators
in the basis associated with\\ $U$-spin reduction}

   First, we consider that the generators $A_{11}$,
$A_{12}$, $A_{21}$, and $A_{22}$ form a basis of the $U$-spin algebra,
which is a compact subalgebra of the noncompact quantum $u_q(2,1)$
algebra, the generators
\beq
U_+=A_{12}\,,\quad U_-=A_{21}\,,\quad U_0=\frac 12
(A_{11}-A_{22})
\eeq
generating the compact $su_q(2)$ subalgebra.

 In the case of $U$-spin reduction, the basis of an unitary irreducible
 representation of the $u_2(2,1)$ algebra can be derived in the same way
 as the basis of $u_q(3)$ algebra  \cite{5}:
\beq
|\langle f_1f_2f_3\rangle m_3UM_U\rangle_q
=\frac{1}{N(k\ell)N(UM_U)}A_{21}^{U-M_U}
P^UA_{23}^kA_{13}^\ell |H\rangle
\eeq
where
\beq
m_3=f_3-k-\ell\,,
\eeq
\beq
U=\frac 12(f_1-f_2-k+\ell)\,,
\eeq
\beq
M_U=\frac 12 (m_1-m_2)\,,\qquad -U\leq M_U\leq U\,,
\eeq
\beq
P^U=\sum_{r=0}^{\infty}(-1)^r\frac{[2U+1]!}{[r]![2U+r+1]!}
A_{21}^rA_{12}^r\qquad
\eeq
is the projection operator for the $su_q(2)$ algebra \cite{17},
\beq
N(UM_U)=\sqrt{\frac{[2U]![U-M_U]!}{[U+M_U]!}},
\eeq
  $N(k\ell)$ are normalization factors, and
\beq
|L\rangle =|\langle f\rangle f_3U_LU_L\rangle , \quad U_L=\frac 12 (f_1-f_2)
\eeq
is the lowest weight vector.
The main distinction between the $u_q(2,1)$ and $u_q(3)$ algebras
lies in the normalization factor $N(k\ell)$. In the Appendix,
it is shown that, in the latter case, the square of the
normalization factor has the form:
%$$
\begin{multline}
N^2(k\ell)=\frac{[k]![\ell]![f_1-f_2-k+\ell+1]![f_1-f_2]!}
{[f_1-f_2-k]![f_1-f_2+\ell+1]!}\\[2ex] %,\times$$
%\beq
\times\,\frac{[f_2-f_3+k-2]![f_1-f_3+\ell-1]!}
{[f_1-f_3-1]![f_2-f_3-2]!}.
\end{multline}
Here, we impose the conventional requirement that the
arguments of all $q$-factorials be nonnegative integers.
This requirement ensures that the square of the norm of basis vectors
is positive. It also follows that a nonzero vector exists only
under the conditions from
$$ f_1\geq f_2, $$
\beq
f_1-f_3\geq 1,
\eeq
\beq
f_2-f_3\geq 2,
\eeq
$$0\leq k\leq f_1-f_2.$$

At the same time, no condition is imposed on the exponent $\ell$
($\ell=0, 1, 2, ...)$, with the result that, in the case of a
$U$-basis, an unitary irreducible representation of the $u_q(2,1)$
algebra is infinite-dimensional.

In \cite{11}, each basis vector of the lowest weight unitary irreducible
representation was characterized by the scheme
\beq
\left|\begin{array}{lll}
m_{13} & m_{23} & m_{33}\\
m_{12}&m_{22}&\\
m_{11}&&\\
\end{array}\right.\Biggr\rangle,
\eeq
where the integers $m_{ij}$ satisfy the conditions:
\beq
m_{13}\geq m_{23}\geq m_{33}\geq 0,
\eeq
\beq
m_{12}\geq m_{13}+1\geq m_{22}\geq m_{23}+1,
\eeq
\beq
m_{12}\geq m_{11}\geq m_{22}.
\eeq
The numbers in the first row in (30) represent the
signature of a unitary irreducible representation of
the $u_q(2,1)$ algebra.
They are related to the components of the lowest weight
by the equations
\beq
f_1=m_{13}+1,
\eeq
\beq
f_2=m_{23}+1,
\eeq
\beq
f_3=m_{33}-2.
\eeq
The numbers in the second row in (30) represent the signature of
a unitary irreducible representation of the $u_q(2)$ subalgebra.
 In our notation,
\beq
m_{12}=f_1+\ell,
\eeq
\beq
m_{22}=f_2+k.
\eeq
The number in the third row is
\beq
m_{11}=U+M_U+m_{22}.
\eeq
 From the condition $f_1\geq f_2$, it follows
\beq
m_{13}\geq m_{23}.
\eeq
The condition $f_2-f_3-2\geq 0$ means that
\beq
m_{23}\geq m_{33}-1.
\eeq
Combining these conditions, we obtain
\beq
m_{13}\geq m_{23}\geq m_{33}-1.
\eeq
At the same time, the condition $0\leq k\leq f_1-f_2$
is equivalent to the constraints
\beq
m_{13}+1\geq m_{22}\geq m_{23}+1.
\eeq
With regard for the allowed values of the exponent $\ell$,
$\ell=0, 1, 2, \ldots$, we derive
\beq
m_{12}\leq m_{13}+1.
\eeq
The condition $-U\leq M_U\leq U$ leads to the constraints
\beq
m_{12}\geq m_{11}\geq m_{22}.
\eeq
A comparison of formulas (42)--(45) with (31)--(33)
demonstrates that our constraints on the structure of basis vectors
 are identical to the constraints  on the values of
$m_{ij}$ in the Gel'fand--Graev schemes, with the only exception that,
in our case, there exists a unitary irreducible representation
for which $m_{23}=m_{33}-1$. This means that there are unitary
irreducible representations corresponding to the Gel'fand--Graev
signature,
\beq
\{m_{13}m_{23}m_{33}\}=\{m_{13}, m_{33}-1,m_{33}\},
\eeq
which are beyond the standard constraints (31). The existence of such
nonstandard discrete series of unitary representations of the $u(2,1)$
algebra was indicated in  \cite{14},  and \cite{15}.
The $u_q(2,1)$ algebra has analogous special series of unitary irreducible
representations.

   Further, it should be noted that at  $f_1=f_2$, in which case $k=0$,
 the condition that the norm $N^2(0\ell)$ is positive requires fulfillment
of inequality
\beq
f_3-f_2-1+\ell>0
\eeq
for all values of $\ell$, including $\ell=1$.
Therefore, the lowest weights corresponding to  $f_1-f_3>0$
are allowed at $f_1=f_2$ (that is, at $m_{13}=m_{33}$).
Therefore, there is an additional nonstandard series of the lowest weight
unitary irreducible representations such that condition (31) is
violated for them. Those are characterized by Gel'fand--Graev signatures
$\{m_{23}-2, m_{23}-2, m_{23}\}$.

   Let us now consider the matrix elements of generators in the $U$ basis.
   In the basis specified by (20), the weight generators $A_{ii}$ ($(i=1,2,3)$)
naturally have a diagonal form are diagonal form, their matrix
elements being given by
\beq
m_1=f_1+\ell-(U-M_U),
\eeq
\beq
m_2=f_2+k+(U-M_U),
\eeq
\beq
m_3=f_3-k-\ell,
\eeq
where
\beq
m_1+m_2+m_3=f_1+f_2+f_3.
\eeq

The action of the generators $A_{12}=U_+$ and $A_{21}=U_-$
are well known from the theory of angular momenta:
\begin{multline}
U_{\pm}|\langle f\rangle m_3UM_U\rangle_q\\[1ex]
=\sqrt{[U\mp M_U][U \pm M_U+1]}\;\,
|\langle f\rangle m_3UM_U\pm 1\rangle_q\,.
\end{multline}
The matrix elements of the generators $A_{13}, A_{23}, A_{31},
A_{32}$ are given by the $q$-analogs of Gel'fand--Graev formulas
\beq
A_{ij} |\langle f\rangle  m_3UM_U\rangle  =\sum_{U'} a_{ij}(m'_3U'M'_U) \;
| \langle f\rangle m'_3U'M'_U\rangle ,
\eeq
where
\beq
U'=U \pm {1/2},\quad M'_U=M_U \pm {1/2}\,
\eeq
and
\beq
a_{ij}(m'_3U'M'_U)={_q\langle} m'_3U'M'_U|A_{ij}|m_3UM_U\rangle_q \,.
\eeq
The list of these matrix elements is given in the Table 1;
their derivation is given in  \cite{18}.
%\newpage

\begin{table} %[p]
\begin{sideways}
%\rotatebox{90}
{\vbox to .9\textwidth{\hsize=.84\textheight
%{\bf Table 1.}\\
\caption{Matrix elements of the generators of the noncompact $u_q(2,1)$
algebra for the unitary irreducible representation
$D^{\{\langle f\rangle +\}}$
associated with the positive discrete series ($U$-basis used here
was derived from the lowest weight vector $|L\rangle$).}
%\vspace{0.5cm}
\extrarowheight=5.2dd
%\noindent
\begin{tabular}{|l|} \hline
$\displaystyle
a_{13}\left(m_3-1,U+\frac12, M_U+\frac12 \right)=
q^{-U+M_U}\left[\frac{[\ell+1][f_1-f_3+\ell][2U+k+2][U+M_U+1]}
{[2U+1][2U+2]}\right]^{1/2}$\\[2ex] \hline
$\displaystyle
a_{23}\left(m_3-1,U+\frac12 , M_U-\frac12 \right)=
%\qquad\qquad
\left[\frac{[\ell+1][f_1-f_3+\ell][2U+k+2][U-M_U+1]}
{[2U+1][2U+2]}\right]^{1/2}$\\[2ex] \hline
$\displaystyle
a_{13}\left(m_3-1,U-\frac 12,M_U+\frac 12\right)=
-q^{U+M_U+1}\left[\frac{[k+1][f_2-f_3+k-1][2U-\ell][U-M_U]}
{[2U][2U+1]}\right]^{1/2}$\\[2ex] \hline
$\displaystyle
a_{23}\left(m_3-1,U-\frac 12,M_U-\frac 12\right)=
\left[\frac{[k+1][f_2-f_3+k-1][2U-\ell][U+M_U]}
{[2U][2U+1]}\right]^{1/2}$\\[2ex] \hline
$\displaystyle
a_{31}\left(m_3+1,U-\frac 12,M_U-\frac 12\right)=-q^{U-M_U}
\left[\frac{[\ell][f_1-f_3+\ell-1][2U+k+1][U+M_U]}
{[2U][2U+1]}\right]^{1/2}$\\[2ex] \hline
$\displaystyle
a_{32}\left(m_3+1,U-\frac 12,M_U+\frac 12\right)=
-\left[\frac{[\ell][f_1-f_3+\ell-1][2U+k+1][U-M_U]}
{[2U][2U+1]}\right]^{1/2}$\\[2ex] \hline
$\displaystyle
a_{31}\left(m_3+1,U+\frac 12,M_U-\frac 12\right)=
q^{-U-M_U-1}
\left[\frac{[k][f_2-f_3+k-2][2U-\ell+1][U-M_U+1]}
{[2U+1][2U+2]}\right]^{\!\!1/2}$\\[2ex] \hline
$\displaystyle
a_{32}\left(m_3+1,U+\frac 12,M_U+\frac 12\right)=
%\qquad\qquad
-\left[\frac{[k][f_2-f_3+k-2][2U-\ell+1][U+M_U+1]}
{[2U+1][2U+2]}\right]^{1/2}$\\[2ex] \hline
\end{tabular}
}}
\end{sideways}
 \end{table}

\section{Basis vectors and matrix elements of the generators
in the basis associated with\\
$T$-spin reduction.}

Let us consider the structure of the unitary irreducible representations
$D^{\langle f\rangle }$ of the lowest weight $(f_1f_2f_3)$ in the case
of the reduction
\beq
u_q(2,1)\rightarrow u_q(1,1)
\eeq
of the $u_q(2,1)$ algebra to the $u_q(1,1)$ subalgebra specified
by generators $A_{22}$, $A_{23}$, $A_{32}$, and $A_{33}$, or to the
$su_q(2)$ subalgebra of a noncompact $T$-spin, the generators in
latter case being
\beq
T_+=A_{23},\quad T_-=A_{32},\quad T_0=\frac 12
(A_{22}-A_{33}).
\eeq
We note that the condition (1) imposed in \cite{11}
on a chain of subalgebras is not satisfied in formulas (56).
For this reason, the results obtained in \cite{11} are not valid
in the case of the $T$-spin basis even for classical $u(2,1)$
algebra, not to mention its deformation $u_q(2,1)$.

 Before proceedings to discuss the $u_q(2,1)$ algebra as a whole,
 it is reasonable to recall general information about the $su_q(1,1)$
 subalgebra and its unitary irreducible representations.
The generators of the $su_q(1,1)$ subalgebra satisfy the well-known
commutation relations
\begin{gather}
[T_0,T_+]=T_+,
\\[1ex] %eeq
%\beq
[T_0,T_-]=-T_-,
\\[1ex] %eeq
%\beq
[T_+,T_-]=[2T_0] .
\end{gather}
Under Hermitian conjugation, they transform as follows:
\beq
T_0^+=T_0,
\eeq
\beq
T_+^+=-T_-.
\eeq
The unitary irreducible representations $D^T$ of the positive discrete
series are infinite-dimensional; the respective $T$-spin is given by
\beq
T=-\frac 12,\; 0,\; \frac 12,\; 1, \; \frac 32,\; 2,\; \ldots
\eeq
The $T$-spin projection $M$ (or the weight of a vector) is an eigenvalue
of the operator of the $T$-spin projection$T_0$, takes the positive values:
\beq
M=T+1,\; T+2,\; \ldots
\eeq
the lowest weight being  $T+1$. We assume that the lowest weight
vector $|H\rangle =|T, T+1\rangle$ is known and that it satisfies the
requirements
\beq
T_-|L\rangle =0,
\eeq
\beq
T_0|L\rangle =(T+1)|L\rangle ,
\eeq
and the normalization condition
\beq
\langle L|L\rangle =1.
\eeq
The basis vectors of a higher weight can be obtained from the lowest
weight vector by the formula
\beq
|TM\rangle =\frac{1}{N(TM)}\;T_+^{M-T-1}\,|TT+1\rangle .
\eeq
The square norm of a vector is derived in this way has a form
$(x=M-T-1)$
\begin{multline}
N^2(TM)=(-1)^x\langle L|A_{32}^xA_{23}^x|L\rangle
=[x][2T+x+1]N^2(T,M+1)\\[1.ex]
= \frac{[-T+M-1]![T+M]!}{[2T+1]!}.
\end{multline}

It can be seen that the condition $N^2(TM)>0$ imposes no
constraints on $x=0,1,2,\ldots$; therefore,
the unitary irreducible representation is infinite-dimensional.
Nonzero matrix elements of the generators in the basis specified
(68) are given by
\beq
\langle TM|T_0|TM\rangle =M,
\eeq
\beq
a_{23}=\langle TM+1|A_{23}|TM\rangle =\{[-T+M][T+M+1]\}^{1/2},
\eeq
\beq
a_{32}= \langle TM-1|A_{32}|TM\rangle =-\{[T+M][-T+M-1]\}^{1/2}.
\eeq
The Casimir operator for the $su_q(1,1)$ algebra has the same form
as for the $su_q(2)$ algebra:
\beq
C_2(su_q(1,1))=T_-T_++[T_0+1/2]^2.
\eeq
All vectors in (64) are the eigenvectors of this operator and
correspond to the same eigenvalue:
\beq
C_2(su_q(1,1))|TM\rangle =[T+1/2]^2|TM\rangle .
\eeq
 We also need the extremal projection operator
 $P^{T}=P^{T}_{T+1.T+1}$ for the discrete series of the lowest
 weight unitary irreducible representation. As in the case of
 the $su_q(2)$ algebra, we seek the expression for the extremal
 projection operator in the form of series:
\beq
P^T=\sum_{r=0} C_rT_+^rT_-^r.
\eeq
In what follows, we apply this projection operator only to those vectors
$|T+1\rangle$ that have a specific weight
 $M=T+1$, but which, in general,  do not have a specific value of $T$-spin
 that is, to vectors that are represented by linear combinations of
%the fo
 \beq
 |T+1\rangle =\sum_{T'}|T',T+1\rangle .
 \eeq

In contrast to the case of $su_q(2)$ algebra, however, the sum over $T'$
is finite in the case under study, because the
inequality $T'\leq M'-1$ must hold for the basis vectors $|T'M'\rangle$
of the positive discrete series. In the case (76) it means
that $T'\leq T$. Hence the variable $T'$ in the sum (76) runs
through the values from $T_{min}=-1/2$ or
 0 up to $T$, depending on whether $T$ is an integer or a half-integer.
 By applying the operator in (75) to the vectors in (76), can show that only
 a finite number of terms in (75) make a non-vanishing contribution, namely,
 those that satisfy $T+1-r\leq 1$ or $1/2$ (that is, $r\leq
T$ or $T+1/2$). Hence, the terms in (75) that involve higher powers $r$
can be disregarded.

 The projection operator $P^T$ satisfies the equations
\beq
T_-P^T=0,
\eeq
\beq
P^T|T,T+1\rangle =|T,T+1\rangle .
\eeq
>From (77) it follows that the coefficients $C_r$ satisfy the recursion relation \beq
C_{r-1}+[r][-2T+r-1]C_r=0.
\eeq
>From this relation, we obtain
\beq
C_r=C_0\frac{[2T-r]!}{[r]![2T]!},\quad r\leq 2T.
\eeq
>From the condition (78), it follows that
\beq
C_0=1.
\eeq
At $r=2T+1$, relation (80) is meaningless, but we have shown above
that we do need the coefficients $C_{r}$ for $r>T$ or $T+1/2$. Thus,
it is sufficient, for our purposes, to use the simple projection operator
\beq
P^T=\sum_{r=0}^{r=2T} \frac{[2T-r]!}{[r]![2T]!}\, T_+^rT_-^r .
\eeq
A projection operator of a more general form can be represented as
\beq
P_{MM'}^T=\frac{(-1)^{-T-M'-1}}{N(TM)N(TM')}\;T_-^{-T+M-1}
P^TT_+^{-T+M'-1}.
\eeq
As a matter of fact, Vilenkin \cite{16} used similar projection
operators (of course, only for $q=1$) long ago
to derive the harmonic projections of polynomials depending on
$n$ Cartesian variables.

  Let us present yet another relation helpful for subsequent computations
\beq
P_{T+1,M}^TP_{M,T+1}^T=\frac{(-1)^{T+M+1}}{N^2(TM)}\; P^TT_-^{M-T-1}
T_+^{M-T-1}P^T=P^T.
\eeq
>From this equation, it follows
\beq
P^TT_-^{M-T-1}T_+^{M-T-1}P^T=(-1)^{M-T-1}N^2(TM)P^T.
\eeq
 We now return to a consideration of the lowest weight unitary
 irreducible representations of the $u_q(2,1)$ algebra. As in the case of
 the $u_q(3)$ algebra, the basis vectors of the unitary
 irreducible representation $D^{\langle f\rangle }$ of the $u_q(2,1)$ algebra
 that correspond to the lowest weight $(f)=(f_{1}f_{2}f_{3})$ will
 be represented in a form
\beq
|\langle f\rangle m_1TM_T\rangle_q=\frac{1}{N(sp)N(TM_T)}\;A_{23}^{M-T-1}P^TA_{13}^s
A_{21}^p|H\rangle ,
\eeq
where
\beq
T=\frac 12(f_2-f_3+p+s-2),
\eeq
\beq
M=T+1,T+2, \dots .
\eeq
The normalization factor $N(TM_T)$ is determined by formula (69),
while the projection operator $P^T$ is given by (82).
The normalization factor $N(sp)$ for the vectors of $T$-spin basis
is calculated by a method similar to that used for the norm of
the vectors of the $U$-spin basis described in the Appendix
(see also \cite{18}). The square of this norm is
\begin{multline}
N^2(sp)=\frac{[s]![p]![f_1-f_2]![f_1-f_3+s-1]!}
{[f_1-f_2-p]![f_1-f_3-1]!}\\[1.5ex]
\times\frac{[f_2-f_3+s-2]![f_2-f_3+p-2]!}{[f_2-f_3-2]!
[f_2-f_3+p+s-2]!}.
\end{multline}
>From the analysis of the norm of the basis vectors, we derive
 the conditions
\beq
f_1\geq f_2,
\eeq
\beq
f_2-f_3-2\geq 0,
\eeq
\beq
0\leq p\leq f_1-f_2 .
\eeq
There is no constraint on the exponent $s$; that is,
$s=0,$ 1, 2, $...$ Since the number of values of the
projections $M$ is infinite, this means that the representations
under study are infinite-dimensional. The constraints
in (91) and (92) on the the signature of unitary irreducible
representations are identical to those obtained for the $U$-spin
basis. For this reason, the classification of the
standard and nonstandard discrete series for the $T$-basis
remains unchanged, as might have been expected.

   Let us now proceed to discuss the matrix elements of
generators. For the weight generators $A_{ii}$
in the $T$-spin basis, only diagonal matrix elements
do not vanish. They are given by
\begin{gather}
a_{11}(m_1TM)=a_{11}(m_{1}m_{2}m_{3})= m_1=f_1-p+s, \\
%\eeq
%\beq
a_{22}(m_{1}TM)=a_{22}(m_{1}m_{2}m_{3})= m_2=f_2+p-T+M-1,\\
%\eeq
%\beq
a_{33}(m_{1} TM)= a_{33}(m_{1}m_{2}m_{3})=m_{3}=f_3-s+T-M+1.
\end{gather}
The matrix elements of the generators $A_{23}=T_+$ and $A_{32}=T_-$
can be determined by formulas (71) and (72).

 The remaining four non-diagonal generators act on the $T$-basis vectors as
follows
\beq
A_{ij}|\langle f\rangle m_{1}TM>_{q}
=\sum_{T'}a_{ij}(m'_{1}T'M')\;|\langle f\rangle m'_{1}T'M'\rangle_{q}
\eeq
where
\beq
T'=T\pm{1/2}, \quad M'=M\pm{1/2},
\eeq
and
\beq
a_{ij}(m'_{1}T'M')={_{q}\langle}m'_{1}T'M'|A_{ij}|m_{1}TM\rangle_q .
\eeq
These matrix elements are presented in Table 2(the derivation
of these expressions is given in \cite{18}.
 %\newpage
\begin{table}[p]
\begin{sideways}
%\rotatebox{90}
{\vbox to .9\textwidth{\hsize=.81\textheight
%{\bf Table 2.}\\
\caption{Matrix elements of the generators of the noncompact $u_q(2,1)$
quantum algebra for the unitary irreducible representation
$D^{\{\langle f\rangle+\}}$ of the positive discrete series ($T$-spin basis
used here was constructed with the aid of the lowest weight vector
$|L\rangle$).}
%\vspace{0.5cm}
%\noindent
\extrarowheight=5.2dd
\begin{tabular}{|l|} \hline
$\displaystyle
a_{12}\left(m_1+1,T+\frac 12,M-\frac 12\right)=
\left[\frac{[s+1][f_1-f_3+s][2T-p+1][-T+M-1]}
{[2T+1][2T+2]}\right]^{1/2}$\\[2ex]  \hline
$\displaystyle
a_{13}\left(m_1+1,T+\frac 12,M+\frac 12\right)=
q^{T-M+1}
\left[\frac{[s+1][f_1-f_3+s][2T-p+1][T+M+1]}
{[2T+1][2T+2]}\right]^{1/2}$\\[2ex]  \hline
$\displaystyle
a_{12}\left(m_1+1,T-\frac 12,M-\frac 12\right)=
\left[\frac{[p][f_1-f_2-p+1][2T-s][T+M]}
{[2T][2T+1]}\right]^{1/2}$\\[2ex] \hline
$\displaystyle
a_{13}\left(m_1+1,T-\frac 12,M+\frac 12\right)=
q^{-T-M}
\left[\frac{[p][f_1-f_2-p+1][2T-s][-T+M]}
{[2T][2T+1]}\right]^{1/2}$\\[2ex]  \hline
$\displaystyle
a_{21}\left(m_1-1,T-\frac 12,M+\frac
12\right)=
\left[\frac{[s][f_1-f_3+s-1][2T-p][-T+M]}
{[2T][2T+1]}\right]^{1/2}$\\[2ex]  \hline
$\displaystyle
a_{31}\left(m_1-1,T-\frac 12,M-\frac 12\right)=
-q^{-T+M-1}
\left[\frac{[s][f_1-f_3+s-1][2T-p][T+M]}
{[2T][2T+1]}\right]^{1/2}$\\[2ex]  \hline
$\displaystyle
a_{21}\left(m_1-1,T+\frac 12,M+\frac 12\right)=
\left[\frac{[p+1][f_1-f_2-p][2T-s+1][T+M+1]}
{[2T+1][2T+2]}\right]^{1/2}$\\[2ex]  \hline
$\displaystyle
a_{31}\left(m_1-1,T+\frac 12,M-\frac 12\right)=
- q^{T+M}
\left[\frac{[p+1][f_1-f_2-p][2T-s+1][-T+M-1]}
{[2T+1][2T+2]}\right]^{1/2}$\\[2ex]  \hline
\end{tabular} }}
\end{sideways}
\end{table}
%\newpage

\section{Weyl coefficients $\langle U|T\rangle_q$ for the positive
discrete series of unitary irreducible representations of the
$u_q(2,1)$ quantum algebra}

By definition, the Weyl coefficient  $\langle U|T\rangle_q$ for an irreducible
 representation $\langle f\rangle $ of the $u_q(2,1)$ quantum algebra
has a form
\begin{multline}
\langle U|T\rangle_q={_q\langle}\{f\}m_3UM_U|\{f\}m_1TM_T\rangle_q\\[1ex]
=\frac{(-1)^{k+\ell}}{N(k\ell)N(UM_U)N(sp)N(TM)}\;
{_q\langle} L|A_{31}^\ell A_{32}^kP^UA_{12}^aA_{23}^bP^T
A_{21}^pA_{13}^s|L\rangle_q , \!\!\!\!\!
\end{multline}
where $|L\rangle$ is the lowest weight vector of the irreducible
representation $\langle f\rangle$:
\beq
a=U-M_U,\qquad b=-T+M-1,
\eeq
and the normalization factors $N(k\ell)$, $N(UM_U)$, $N(sp)$,
and $N(TM)$ and the projection operators $P^U$ and $P^T$ were
defined in the foregoing.
Since the weight in the left hand side of Eq.~(99 ) for the
matrix element is equal to the weight of the right hand side,
we conclude that the parameters $k$ and $\ell$ are related to $s$
and $p$ by the equations
\beq
U-M_U=p-s+\ell,
\eeq
\beq
T-M+1=s-\ell-k.
\eeq

   The computation of the above matrix element is performed by
making use of the commutation relations between the generators raised
to a power. The scheme of the computations is identical to that
in the case of the $u_q(3)$ algebra \cite{19}.
Taking into account the explicit form of projection operator
$P^U$, we arrive at
\begin{multline}
B={_q\langle}L|A_{13}^\ell A_{32}^k P^U
A_{21}^rA_{12}^{a}A_{23}^bP^TA_{21}^pA_{13}^s|L\rangle_q \\ %=\,
=\sum_r(-1)^r\frac{[2U+1]!}{[r]![2U+r+1]!}\;
%\times
{_q\langle} L|A^{\ell}_{31}A^{k}_{32}A^{r}_{21}A^{a+r}_{12}A^b_{23}
P^TA^p_{21}A^s_{13}|L\rangle_q    .
\end{multline}
With the aid of commutation relations, we transfer the operator
$A^r_{21}$ in the matrix element to the left until it appears
immediately after vector $\langle L|$ and consider that $\langle L|A_{21} =0$.
The expression for $B$ then takes the form
\beq
B=\sum_r\frac{[2U+1]![k]!}{[r]![k-r]![2U+r+1]!}\:B_1,
\eeq
where
\beq
B_1={_q\langle} L|A_{31}^{\ell+r}A_{12}^{r+a}A_{32}^{k-r}A_{23}^b
P^TA_{21}^{p}A_{13}^{s}|L\rangle_q .
\eeq
To compute the matrix element $B_1$, the generators $A_{32}$
must be transferred to the right until they appear immediately
before the projection operator $P^T$, whereupon the equation
$A_{32}P^T=0$ is taken into account. As a result the matrix
element $B_1$ reduces to the expression
 \beq
 B_1= \frac{[k]![b]!}{[k-r]![b-k+r]!}\;\prod_{t}
 [f_3-f_2-a+b-k-\ell-r-t]\;B_2,
 \eeq
 where
 \beq
 B_2={_q\langle}
L|A^{\ell+r}_{31}A^{a+r}_{12}A^{b-k+r}_{23}P^TA^p_{21}A^s_{13}
    |L\rangle_q.
 \eeq
 Considering that, in the case of $T$-spin basis,
 \beq
 2T=f_2-f_3+p+s-2\,
 \eeq
 and that all factors in the product $\prod\limits_t$ are negative, we recast
 the product into the form
\beq
\prod_{t}[f_3-f_2-a+b-k-\ell-r-t]=
(-1)^{k-r}\frac{[f_2-f_3+p+\ell+k-1]!}{[f_2-f_3+p+l+r-1]!}.
\eeq
Further, we transfer the generators $A^{b-k+r}_{23}$ in the expression
for the matrix element $B_2$ to the left until they appear immediately
after the vector $\langle L|$, which annihilates them, and consider that
$\langle L|A^y_{21}=\delta_{y,0}\langle L|$. As a result we have
\beq
B_2=\frac{[a+r]!}{[a-b+k]!}\;{_q\langle} L|A^{\ell+r}_{31}A^{b-k+r}_{13}
A^{a-b+k}_{12}P^TA^p_{21}A^s_{13}|L\rangle_q .
\eeq
The commutation of generators $A^{\ell+r}_{31}$ and $A^{b-k+r}_{13}$
whereupon the condition
$\langle L|A^z_{13}=\delta_{z,0}\langle L|$ is taken into
account, gives the ultimate expression for the matrix element $B_2$:
\begin{multline}
 B_2=(-1)^{\ell+r+s}\:\frac{[a+r]![\ell+r]![f_1-f_3+\ell+r-1]!}
 {[p]![s]![f_1-f_3+s-1]!}  \\[1.5ex] %times\,
 \times
{_q\langle} L|A^s_{31}A^p_{12}P^TA^p_{21}A^s_{13}|L\rangle_q \\[1.5ex]
 =(-1)^{\ell+r}\:\frac{[a+r]![\ell+r]![f_1-f_3+\ell+r-1]!}
 {[p]![s]![f_1-f_3+s-1]!}\:N^2(sp).
 \end{multline}
Combining the above results, we reduce the expression for the
Weyl coefficients $\langle U|T\rangle_q$ of the form
\begin{multline}
\langle U|T\rangle_q
=\left\{\frac{[2U+1][2T+1][k]![-T-1+M][U+M_U]![T+M]!}
{[p]![s]![\ell]![U-M_U]![f_1-f_3+s-1]!}\right. \\[1.5ex]
%\times,\,
\left.\times\frac{[f_1-f_2-k]![f_1-f_2+\ell+1]![f_2-f_3+s-2]!
[f_2-f_3+p-2]!}{[f_1-f_2-p]![f_2-f_3+k-2]![f_1-f_3+\ell-1]!}\right\}^{1/2}
                    \\[1.5ex]
%\times\,
\times\sum_r(-1)^r\frac{[U-M_U+r]![\ell+r]![f_1-f_3+\ell+r-1]!}
 {[r]![2U+r+1]![k-r]![\ell-s+r]![f_2-f_3+p+\ell+r-1]!} .
\end{multline}
The substitution $r=k-n$ allows to rewrite the last formula as follows
\begin{multline}
\langle U|T\rangle_q =
\left\{\frac{[2U+1][2T+1][k]![-T+M-1]![U+M_U]![T+M]!}
{[s]![p]![\ell]![U-M_U]![f_1-f_3+s-1]!}
    \right. \\[1.1ex] %times\
\times \left.
\frac{[f_1-f_2-k]![f_1-f_2+\ell+1]![f_2-f_3+s-2]![f_2-f_3+p-2]!}
{[f_1-f_2-p]![f_2-f_3+k-2]![f_1-f_3+\ell-1]!}\right\}^{1/2}
  \\[1.1ex]   %\times\
\times\sum_n \frac{(-1)^{k+n} [U-M_U+k-n]![\ell+k-n]!}
{[n]![k-n]![2U+1+k-n]![\ell-s+k-n]!}
   \\[1.1ex] \times \frac{[f_1-f_3+\ell+k-n-1]!}{[f_2-f_3+p+\ell+k-n-1]!} .
\end{multline}

\section{Relation between the $q$-Weyl coefficients for the $u_q(2,1)$
quantum algebra and $q$-Racah coefficients for the $su_q(2)$
quantum algebra.}

The explicit expression for the $q$-Weyl coefficient for the $u_q(3)$
algebra was obtained in \cite{19}, and its relation to the Racah
coefficient for the $su_q(2)$ quantum algebra was established there.
  Here, we show that the expression (113)for the $q$-Weyl coefficient
for the $u_q(2,1)$ quantum algebra can also be related to the $q$-Racah
coefficients for the $su_q(2)$ quantum algebra.
Our consideration is based on one of five general formulas in \cite{19}
for the $q$-Racah coefficient for the $su_q(2)$ quantum algebra [namely,
formula (5.31) in \cite{19}]:
\begin{multline}
U_q(abed;cf)=
(-1)^{a+d-c-f}\left\{\frac{[2c+1][2f+1][a+b+c+1]!}
{[a+e+f+1]![c+d+e+1]!}\right. \\[1.5ex]
\times\frac{[b+d+f+1]![a-b+c]!
[-a+b+c]![a+e-f]![b-d+f]!}{[a+b-c]!
  [a-e+f]![b+d-f]![c+d-e]![c-d+e]!}\\[1.5ex]
\left.\times\frac{[-b+d+f]![-c+d+e]!}{[-a+e+f]!}\right\}^{1/2} \\[1.5ex]
%\times
\times\sum_{n}\frac{(-1)^{n}[2b-n]![b+c-e+f-n]!}
{[n]![-a+b+c-n]![b-d+f-n][a+b+c+1-n]!}\\[1.5ex]
\times\frac{1}{[b+d+f+1-n]!} .
\end{multline}

A comparison of the expression (114) with formula (113)
for the Weyl coefficient for the positive discrete series
of the representations of $u_q(2,1)$ reveals the summands in
the two formulas coincide, provided
$$
a=T=\frac12(f_2 -f_3 +p+s-2),\;\;\;  b=j_3=\frac12 (\ell+k) ,
$$
\beq
c=j_2= \frac12(f_2 -f_3+p-s+\ell+k-2)), %\;\;\;\;.
\eeq
$$
d=U=\frac12(f_1-f_2+\ell-k),\;\;\;\;\;\;
e =j_1=\frac12(f_1 -f_3-p+s-2),
$$
\beq
 f=j=\frac12 (f_1 -f_2).
\eeq
Further, the substitution of parameters  $a,b,c,d,e$ and $f$
from the formulas (115) and (116) gives the relation between the $q$-Weyl
coefficient (113) for the $u_q(2,1)$ quantum algebra and $q$-Racah
coefficient for the $su_q(2)$ quantum algebra,
\begin{multline}
\langle U\mid T\rangle_q
=(-1)^{s}\;{\sqrt{\frac{[2U+1][2T+1]}{[2j_2+1][2j+1]}}\;U(Tj_3j_1U;j_2j)_q }\\
%$$
%\beq
=(-1)^{k}\, U(j_1j_2jj_3;UT)_{q}.
\end{multline}
% \newpage

\section{Conclusion}

In this study, the projection operators for the $su_q(2)$
subalgebra have used to explore the positive discrete series
of unitary irreducible representations of the noncompact $u_q(2,1)$
quantum algebra. The $q$-analog of the Gel'fand--Graev formulas
has been derived in the basis associated with the reduction
$u_q(2,1)\rightarrow su_q(2)\times u(1)$. It seems that the reduction
$u_q(2,1)\rightarrow u(1)\times su_q(1,1)$ for the discrete series of
the lowest weight representations has been considered for the first time
in the present study. With the aid of the projection operator for the
$su_q(1,1)$ subalgebra, we constructed the basis of the representation for
this reduction and calculated the matrix elements of the generators.
We have obtained analytic expressions for the elements of
the transformation brackets
$\langle U|T\rangle_q$ relating the $U$-spin and $T$-spin
bases of the lowest weight irreducible representations. By the
analogy with $q$-Weyl coefficients for the $u_q(3)$ algebra
\cite{19}, they can be called the $q$-Weyl coefficients for the
noncompact $u_q(2,1)$ algebra. It has been explicitly shown that
these $q$-Weyl coefficients are equivalent (apart from phase factor )
to specific $q$-Racah coefficient for the $u_q(2)$ algebra or are
proportional to the $q$--6j symbol for the $su_q(2)$ algebra.
The negative discrete series was discussed by us in \cite{20}.
The intermediate discrete series requires a dedicated investigation,
and this will be done in our further publication.

\vspace{0.5cm}
This work was supported by the Russian Foundation of Basic Research
(project No 02-01-00668).

%\vspace{2cm}
%\newpage
%{\bf Appendix.}
%\begin{center}
%\underline
\chapappendix{Normalization of the $U$-spin  basis vectors
of the $u_q(2,1)$ algebra  %\\
(positive discrete series).}
%\end{center}

The structure of the $U$-basis vectors is described by formulas
(19)--(26).

Here, we use the transformation properties of the ``noncompact''
generators under Hermitian conjugation and the properties of
projection operator $P^U$:
$$(P^U)^+=P^U,\eqno{({\mbox A}.1)}$$
$$(P^U)^2=P^U,\eqno{({\mbox A}.2)}$$
$$A_{12}^{U-M_U}A_{21}^{U-M_U}P^U=N^2(UM_U)P^U.\eqno{({\mbox
A}.3)}$$
With allowance for these formulas, the square of the norm $N^2(k\ell)$
takes a form:
$$N^2(k\ell)=(-1)^{k+\ell}\,\langle L|A_{31}^\ell\tilde
A_{32}^kP^U A_{23}^kA_{13}^\ell|L\rangle ,\eqno{({\mbox A}.4)}$$
where
$$\tilde A_{31}=A_{32}A_{23}-q A_{23}A_{32}.
\eqno{({\mbox A}.5)} $$
Since, by definition, the relation
$$A_{31}=A_{32}A_{21}-q^{-1}A_{21}A_{32} \eqno{({\mbox A}.6)}$$
holds we can represent the generator $\tilde A_{31}$ in the form
$$\tilde A_{31}=A_{31}-(q-q^{-1})A_{21}A_{32}.
\eqno{({\mbox A}.7)}$$
>From the relations
$$A_{12}P^U=P^UA_{21}=0 \eqno{({\mbox A}.8)}$$
it follows that
$$N^2(k\ell)=(-1)^{k+\ell}\,\langle L|A_{23}^\ell A_{13}^kP^UA_{31}^k
A_{32}^\ell|L\rangle .\eqno{({\mbox A}.9)}$$
A straightforward computation of $N^2(k\ell)$ by transferring of lowering
generators to the lowest vector $|L\rangle$ is rather cumbersome.
In view of this, we will try to construct a recursion relation
between the expressions for $N^2(k \ell)$ for various values of
$k$ and $\ell$, bearing in mind that
$$\langle L|P^U|L\rangle =\langle L|L\rangle =1.\eqno{({\mbox A}.10)}$$
We begin by establishing a relation between $N^2(k\ell)$ and $N^2(k-1,\ell)$.
In the expression for $N^2(k\ell)$, we replace, for this purpose,
$A_{32}^k$ by $A_{32}^{k-1}P^{U+1/2}A_{32}$.
This is legitimate because, in (A.9), the projection operator $P^{U+1/2}$
taken in this combination is equivalent to the identity operator.
Indeed, we have
$$P^{U+1/2}A_{32}P^U=\sum_{r}(-1)^r\,\frac{[2U+2]!}
{[r]![2U+r+2]!}\:A_{21}^rA_{12}^rA_{32}P^U %\,=$$
%$$
=A_{32}P^U,\eqno{({\mbox A}.11)}$$
since the generators $A_{12}^r$ and $A_{32}$ commute  and since
$A_{12}^r P^U=\delta_{r,0}P^U$.
We now consider the application of the generator $A_{32}$ to
the projection operator $P^U$:
\begin{multline}
A_{32}P^U=\sum_r(-1)^r
\frac{[2U+1]!}{[r]![2U+r+1]!}\:A_{32}A_{21}^{r}A_{12}^r=\\[1.ex] %\,$$
                                %$$
=\sum_r (-1)^r \frac{[2U+1]!}{[r]![2U+r+1]!}\,(q^{-r}A^r_{21}A_{32}
+[r]A_{21}^{r-1}A_{23})\,A^r_{12}. \tag{A.12}%\eqno{({\mbox A}.12)}
\end{multline}
>From here on, we use the commutation relations from \cite{1,19}
for generators raised to a power. In view of the relation
$P^{U+{1/2}}A_{21}=({A_{12}}P^{U+{1/2}})^{+}=0$, the application
of this operator on the projection operator $P^{U+{1/2}}$ from the left
yields
$$
P^{U+1/2}A_{32}P^U=P^{U+1/2}\left(A_{32}-\frac{A_{31}A_{12}}
{[2U+2]}\right),
\eqno{({\mbox A}.13)}
$$
where
$$[2U+2]=[f_1-f_2-k+\ell+2]. \eqno{({\mbox A}.14)}$$
As a result, the square of the norm becomes
$$
N^2(k\ell)=(-1)^{k+\ell}\left\langle\! L\left|A_{31}^\ell A_{32}^{k-1}
P^{U+1/2} %\,\times$$
%$$\times\,
\left(A_{32}-\frac{A_{32}A_{21}}
{[f_1-f_2-k+\ell+2]}\right)A_{23}^kA_{13}^\ell\right|
L\!\right\rangle_{\vphantom{{q_q}}} \!.
\eqno{({\mbox A}.15)}$$
Commuting the generators $A_{32}$ and $A_{23}^k$, we arrive at
$$
A_{32}A_{23}^kA^{\ell}_{13}\,|L>=[k][f_3-f_2-k-\ell+1]A^{k-1}_{23}A^{\ell}_{13}\,|L>.
\eqno{(\mbox{A}.16)}
$$
and
$$A_{31}A_{12}A_{23}^kA_{13}^{\ell}\,|L\rangle =[k]A_{31}A^{k-1}_{23}
A^{\ell+1}_{13}\,|L>.\eqno{({\mbox A}.17)}
$$

The commutation of the generators $A_{32}$ and $A^{\ell}_{31}$
makes it possible to derive the relation
$$
A_{32}A_{13}^{\ell}=\left(A_{13}^{\ell}-[\ell]q^{\ell-1}A_{13}^{\ell-1}
A_{12}q^{-(A_{22}-A_{33}-1)}\right).\eqno{({\mbox A}.18)}
$$
Transferring the generator $A_{31}$ to the right until it
appears immediately in front of the lowest weight $|L\rangle $,
which annihilates it, we obtain
$$P^{U+{1/2}}A_{32}A_{12}A_{23}^{k}A_{13}^{\ell}|L\rangle =
\,[k][\ell+1][f_3-f_1-\ell]P^{U+{1/2}}A_{23}^{k-1}A_{13}^{\ell}|L\rangle ;
\eqno{({\mbox A}.18)}$$
therefore, we have
%$$
\begin{multline}
P^{U+{1/2}}A_{32}P^UA_{23}^{k}A_{13}^{\ell}\,|L\rangle
\\[1.5ex]
=[k]([f_3-f_2-k-\ell+1]-\frac{[\ell+1][f_3-f_1-\ell]}{[f_1-f_2-k+\ell+2]}\,
P^{U+{1/2}}A^{k-1}_{23}A^{\ell}_{13}\,|L\rangle  \\[1.5ex] %=
=-\frac{[k][f_1-f_2-k+1][f_2-f_3+k-2]}{[f_1-f_2-k+\ell+2]}\:P^{U+{1/2}}
A^{k-1}_{23}A^{\ell}_{13}\,|L\rangle .
%%%%\eqno{({\mbox A}.19)}
\tag{A.19}
\end{multline}
%$$
Thus, the square of the norm, $N^2(k\ell)$, takes the form
\begin{multline}
N^2(k\ell)=(-1)^{k+\ell-1}\frac{[k][f_1-f_2-k+1]
[f_2-f_3+k-2]}{[f_1-f_2-k+\ell+2]} \\[1.5ex] %\,\times
%$$
%$$
\times
\langle L|A_{31}^\ell A_{32}^{k-1}P^{U+1/2}A_{23}^{k-1}
A_{23}^\ell|L\rangle . \tag{A.20}
%\eqno{(\mbox{A}.20)}
\end{multline}
In other words we have derived a recursion relation between
$N^2(k\ell)$ and $N^2(k-1,\ell)$,
$$N^2(k\ell)=\frac{[k][f_1-f_2-k+1][f_2-f_3+k-2]}
{[f_1-f_2-k+\ell+2]}\,N^2(k-1,\ell).\eqno{({\mbox A}.21)}$$
The recursion relation
$$N^2(0\ell)=(-1)^\ell\,\langle H|A_{31}^\ell P^UA_{13}^\ell|H\rangle %\,=$$
 %$$
=[\ell][f_1-f_3+\ell-1]\,N^2(0,\ell-1)
\eqno{({\mbox A}.22)}$$
can be obtained in a similar way.

Using these recursion relations, we arrive at an ultimate expression for
for the square of the norm in (A.9):
$$N^2(k\ell)\,=\frac{[k]![\ell]![f_1-f_2]![f_1-f_2-k+\ell+1]!
[f_2-f_3+k-2]![f_1-f_3+\ell-1]!}
{[f_1-f_2-k]![f_1-f_2+\ell+1]![f_1-f_3-1]![f_2-f_3-2]!}
.\eqno{(\mbox{A}.23)}
$$
% \newpage

%\normallatexbib
%\begin{thebibliography}{99}

\end{document}